\documentclass[a4paper,11pt]{amsart}

\usepackage{amsfonts}
\usepackage{latexsym}
\usepackage{amssymb}
\usepackage[all]{xy}

\newtheorem{thm}{\sc theorem}
\newtheorem{dfn}[thm]{\sc definition}
\newtheorem{lem}[thm]{\sc lemma}
\newtheorem{lemma}[thm]{\sc lemma}

\newtheorem{cor}[thm]{\sc corollary}

\theoremstyle{definition}
\newtheorem{ex}[thm]{\sc Example}
\newcommand{\D}{\mathcal{D}}

\newcommand{\lotimes}{\otimes^\mathrm{L}}

\newcommand{\Cbar}{C^\mathrm{bar}}

\newcommand{\Aop}{{A^\mathrm{op}}}

\newcommand{\rmod}{R\mbox{-}\mathbf{Mod}}
\newcommand{\smod}{S\mbox{-}\mathbf{Mod}}
\newcommand{\amod}{A\mbox{-}\mathbf{Mod}}
\newcommand{\amoda}
{A^e\mbox{-}\mathbf{Mod}}

\newcommand{\umod}{U\mbox{-}\mathbf{Mod}}

\newcommand{\modu}{U^\mathrm{op}\mbox{-}\mathbf{Mod}}

\newcommand{\moda}{A^\mathrm{op}\mbox{-}\mathbf{Mod}}
\newcommand{\lact}{{\,{\scriptstyle\rhd}\,}}
\newcommand{\ract}{{\,{\scriptstyle\lhd}\,}}
\newcommand{\blact}{{\,{\scriptstyle\blacktriangleright}\,}}

\newcommand{\bract}{{\,{\scriptstyle\blacktriangleleft}\,}}

\newenvironment{pf}
{{{\bf Proof.}}}{\hfill$\Box$\\[1mm]}

\begin{document}

\title{Duality and products in 
algebraic (co)homology
theories}
\author{Niels Kowalzig} 
\author{Ulrich Kr\"ahmer} 

\address{N.K.: Utrecht University, 
Department of Mathematics,
P.O. Box 80.010,
3508TA Utrecht,
The Netherlands}

\email{nkowalzi@science.uva.nl}

\address{U.K.: University of Glasgow,
Department of Mathematics, University 
Gardens, Glasgow G12 8QW, Scotland}

\email{ukraehmer@maths.gla.ac.uk}

\begin{abstract}
The origin and interplay of 
products and dualities 
in algebraic (co)homology theories 
is ascribed to a $\times_A$-Hopf algebra
structure on the relevant 
universal enveloping algebra. This
provides a unified treatment for example 
of results by Van den
Bergh about Hochschild (co)homology and by 
Huebschmann about Lie-Rinehart
(co)homology. 
\end{abstract}
\maketitle

\section{Introduction}
Most classical 
(co)homology theories of algebraic 
objects such as groups or Lie, 
Lie-Rinehart or associative algebras
can be realised as
\begin{equation}\label{anfang}
		  H^\bullet(X,M):=
		  \mathrm{Ext}^\bullet_U(A,M),\quad
		  H_\bullet(X,N):=
		  \mathrm{Tor}^U_\bullet(N,A)
\end{equation} 
for an augmented ring $X=(U,A)$ 
(a ring with a distinguished left module)
that is functorially attached to a
given object.
The cohomology coefficients 
are left $U$-modules $M$ and
those in homology  
are right $U$-modules $N$.

Our aim here is to clarify the origin 
and interplay of multiplicative structures 
and dualities between such
(co)homology groups, and to 
provide 
a unified treatment of results  
by Van den Bergh on Hochschild
(co)homology \cite{VdB:ARBHHACFGR} 
and by Huebschmann on 
Lie-Rinehart (co)homology 
\cite{Hue:DFLRAATMC}. 
The key 
concept involved is that of a $\times_A$-Hopf 
algebra introduced by 
Schauenburg \cite{Schau:DADOQGHA}.

The main results can be summarised as
follows:
\begin{thm}\label{main}
For any $A$-biprojective $\times_A$-Hopf
 algebra $U$ there is a functor
$$
		\otimes : \umod \times \modu
		  \rightarrow \modu  
$$
that induces for
$M \in \umod,N \in \modu$
and $m,n \ge 0$
natural products
$$
		  \smallfrown \,\, : 
		  \mathrm{Ext}^m_U(A,M) \times
		  \mathrm{Tor}^U_n(N,A)
		  \rightarrow \mathrm{Tor}^U_{n-m}
		  (M \otimes N,A). 
$$

If $A \in \umod$ admits a finitely
generated projective resolution of
finite length and
there exists $d \ge 0$ with
$\mathrm{Ext}_U^m(A,U)=0$  
for $m \neq d$, then 
there is a canonical element
$$
		  [\omega] \in
 \mathrm{Tor}^U_d(A^*,A),\quad
A^*:=\mathrm{Ext}_U^d(A,U)$$
such that
for $m \ge 0$
and $M \in \umod$ with
$\mathrm{Tor}^A_q(M,A^*)=0$ for $q>0$
$$
		  \cdot \smallfrown [\omega] :
		  \mathrm{Ext}_U^m(A,M)
 \rightarrow 
		  \mathrm{Tor}^U_{d-m}
		  (M \otimes A^*,A)
$$
is an isomorphism.
\end{thm}

As we will recall below,  
$\times_A$-bialgebras
and $\times_A$-Hopf algebras 
generalise bialgebras and 
Hopf algebras towards 
possibly noncommutative base algebras $A$. 
Besides Hopf algebras, both the
universal enveloping algebra $U(A,L)$
of a Lie-Rinehart algebra $(A,L)$ and the
enveloping algebra $A^e=A \otimes_k
A^\mathrm{op}$ of an associative algebra
$A$ are $\times_A$-Hopf algebras, see
Section~\ref{secex}.

For any $\times_A$-bialgebra $U$, the
base algebra $A$ carries a left $U$-action
and the category $\umod$ of left $U$-modules 
is monoidal 
with unit object $A$. But only for 
$\times_A$-Hopf algebras one has
a canonical operation 
$\otimes$ as in Theorem~\ref{main}
which turns $\modu$ into a module
category over 
$(\umod,\otimes,A)$ 
(Lemma~\ref{RechtsLinks}).

Any $\times_A$-Hopf algebra carries two
left and two right actions of the base
algebra, all commuting with each other. 
The biprojectivity assumed in
Theorem~\ref{main} refers to the
projectivity of two of these, see
Section~\ref{soco}. Under this condition, we 
can use the elegant formalism of suspended
monoidal categories from
\cite{Sua:THHAFTACOCP} 
to define for $M,N \in \umod$ and
$P \in \modu$ products
$$
 \smallsmile \,\,: 
			H^m(X,M) \times 
		  H^n(X,N)  \rightarrow 
		  H^{m+n}(X,M \otimes
		  N),
$$ 
$$
 \smallfrown \,\,: H^n(X,N) \times 
		  H_p(X,P) \rightarrow
		 H_{p-n} (X,N \otimes P),
$$
where we again use the abbreviations
from (\ref{anfang}) 
(cf.~Sections~\ref{absch1} and~\ref{absch2}).

In the last part of Theorem~\ref{main},
$A^*=H^d(X,U)=
\mathrm{Ext}^d_U(A,U)$ is a right 
$U$-module via right multiplication 
in $U$, and if we define the functor
$$
		  \hat{} \,\,: \umod \rightarrow
		  \modu,\quad
		  M \mapsto \hat M:=M \otimes A^*,
$$
then the statement can be rewritten
as an isomorphism    
$$
		H^m(X,M) \simeq 
		H_{\mathrm{dim}(X)-m}(X,\hat M),\quad
		\mathrm{dim}
		(X):=\mathrm{proj.dim}_U(A) 
$$
that is given as in topology by the cap
product with the fundamental class
$[\omega] \in H_{\mathrm{dim}(X)}(\hat A)$ which
corresponds under the duality to  
$\mathrm{id}_A \in 
H^0(A)=\mathrm{Hom}_U(A,A)$. 
For $M=A$ this simply means that the
$H^\bullet(A)$-module 
$H_\bullet(A^*)$ is
free with generator $[\omega]$.

Theorem~\ref{main} is well-known 
in group and Lie
algebra (co)homology. For $U=A \otimes_k
A^\mathrm{op}$ it reduces to Van den
Bergh's result \cite{VdB:ARBHHACFGR} that stimulated
a lot of recent research, see e.g.~
\cite{BroZha:DCATHCFNHA,Dol:TVDBDATMSOAPV,Gin:CYA,
LauRic:TPDFSQPA}. Note that
we do not need 
Van den Bergh's invertibility 
assumption about $A^*$, which says 
that $\hat {}$ is an
equivalence. However, it is satisfied
for many well-behaved algebras [ibid.]
and implies the condition 
$\mathrm{Tor}^A_q(M,A^*)=0$ for
arbitrary $A$-bimodules $M$ (since
invertible bimodules are finitely
generated projective as one-sided
modules from either side). 
For Lie-Rinehart
algebras $(A,L)$, Theorem~\ref{main} is due 
to Huebschmann 
\cite{Hue:DFLRAATMC}, and we find   
the general setting helpful for
example to understand the different 
roles of left and right modules that
he has observed. As
Huebschmann has showed, the
conditions of Theorem~\ref{main} are
satisfied whenever $L$ is
finitely generated projective over $A$, and $A^*$ coincides as
an $A$-module with $\Lambda^d_AL$ and is
in particular projective, so also here we
have $\mathrm{Tor}^A_q(M,A^*)=0$ for
arbitrary $(A,L)$-modules $M$.

We were also motivated by 
the current discussion 
of the numerous bialgebroid
generalisations of Hopf algebras, see
\cite{Boe:HA}. Several authors have 
raised the question 
where Lie-Rinehart algebras fit in. They
were shown in
\cite{Xu:QG,MoeMrc:OTUEAOALRA}
to be $\times_A$-bialgebras, see also 
\cite{KhaRan:CCOEHA,Hue:TUHAAWAHLRA}.
Here we remark that they are
in fact always $\times_A$-Hopf algebras,
but not necessarily  
Hopf algebroids in the sense of B\"ohm
and Szlach\'anyi (Example~\ref{bumpy},
this answers a question of B\"ohm 
\cite{Boe:HA}).
So both these examples and the
applications in homological algebra 
clearly demonstrate the relevance of 
the intermediate concept of a 
$\times_A$-Hopf algebra.

Theorem~\ref{main} could be generalised
to differentially graded $\times_A$-Hopf
algebras, sheaves of such, or suitable
abstract monoidal
categories. One can also drop the
condition $\mathrm{Ext}^n_U(A,U)=0$ 
for $n \neq d$ and the assumption that
$\mathrm{Tor}^A_q(M,A^*)=0$. Then one obtains
for a bounded below chain complex $M$
over $\umod$ an isomorphism
$\mathrm{RHom}_U(A,M) \simeq
(M \otimes_A^L \mathrm{RHom}_U(A,U))
\otimes_U^L A$.\\

N.K.~is supported by the NWO
through the GQT cluster. U.K.~is supported by the EPSRC 
fellowship EP/E/043267/1 and partially by the 
Polish Government Grant N201 1770 33. We
thank Andy Baker, Gabriella B\"ohm, Ken
Brown, Henning Krause and 
Valery Lunts for discussions.  

\section{Preliminaries on 
$\times_A$-Hopf algebras}
\subsection{Some conventions}\label{soco}
Throughout this paper, ``ring'' means 
``unital and associative ring'', and we fix   
a commutative ring $k$.
All other algebras, modules etc.~will 
have an underlying structure of a
$k$-module.

Secondly, we fix a
$k$-algebra $A$, i.e., a ring with a 
ring homomorphism 
$ \eta_A : k \rightarrow Z(A)$ to  
its centre.
We denote by $\amod$ the category of 
left $A$-modules, 
by $A^\mathrm{op}$ the 
opposite and by
$A^e:=A \otimes_k A^\mathrm{op}$ 
the enveloping algebra 
of $A$. Thus left 
$A^e$-modules are $A$-bimodules
with symmetric action of $k$.

Our main object is finally an algebra $U$
over $A^e$, where we now refer to the
less standard notion of 
an algebra over a possibly noncommutative base
algebra:  
$U$ is a $k$-algebra with a $k$-algebra 
homomorphism $ \eta=\eta_U : A^e \rightarrow
U$. This gives rise to a
forgetful functor 
$\umod \rightarrow \amoda$ using which
we consider every 
$U$-module $M$ also as an $A$-bimodule
with actions
\begin{equation}\label{brot}
		  a \lact m \ract b := \eta (a \otimes_k b)m,\quad
		  a,b \in A,m \in M.
\end{equation} 
Similarly, every right $U$-module $N$ is also 
an $A$-bimodule via
\begin{equation}\label{salz}
		  a \blact m \bract b := 
			n \eta (b \otimes_k a),\quad
		  a,b \in A,n \in N.
\end{equation} 
In particular, $U$ itself carries 
two left and two right 
$A$-actions all commuting 
with each other.
Usually
we consider $U$ as an $A^e$-module using 
$a \lact u \ract b$, and otherwise
we write e.g.~${}_\blact U_{\ract}$ to denote 
which actions are considered. Since this
will be  repeatedly a necessary
technical condition, we define:
\begin{dfn}
For an $A^e$-algebra $U$ we call 
$M \in \umod$
$A$-biprojective if both ${}_\lact M \in \amod$ 
and $M_\ract \in \moda$ are 
projective modules.
\end{dfn} 
 
\subsection{$\times_A$-bialgebras \cite{Tak:GOAOAA}}  
Consider an 
$A^e$-algebra $U$ as above which 
is also a coalgebra 
in the monoidal category 
$\amoda$. That is, there are maps
$$
		  \Delta : U \rightarrow U \otimes_A U,\quad
		  \varepsilon : U \rightarrow A
$$
satisfying the usual coalgebra 
axioms (see e.g.~\cite{Boe:HA} 
for the details), where
\begin{equation}\label{daoben}
		  U \otimes_A U=U \otimes_k U/
		  \mathrm{span}_k
		  \{u \ract a \otimes_k v-u \otimes_k 
		  a \lact v\,|\,
		  a \in A,u,v \in U\}.
\end{equation}

For $A=k$ one calls $U$ a bialgebra if 
$\Delta$ and $\varepsilon$ are algebra
homomorphisms, but in general there is no 
natural algebra structure on
$U \otimes_A U$. The way out of this problem
was found by Takeuchi \cite{Tak:GOAOAA}
and involves the embedding
\begin{equation}\label{iota}
		  \iota : U \times_A U \rightarrow 
		  U \otimes_A U,
\end{equation} 
where $U \times_A U$ is the centre of
the 
$A$-bimodule ${}_\blact  U_\ract
\otimes_A {}_\lact U_\bract$:
$$
		  U \times_A U:=
		  \Bigl\{\sum_i u_i \otimes_A v_i \in 
		  U \otimes_A U\,|\,
		  \sum_i a \blact u_i \otimes_A v_i=
		  \sum_i u_i \otimes_A v_i \bract a
		  \Bigr\}.
$$
The product of $U$ turns this into 
an algebra over $A^e$, with 
$$
		  \eta_{U \times_A U} : A^e \rightarrow 
		  U \times_A U,\quad
		  a \otimes_k b \mapsto 
		  \eta (a \otimes_k 1) \otimes_A 
		  \eta(1 \otimes_k b).
$$ 
Similarly, 
$A$ is an algebra over $k$, but 
not over $A^e$ in general. 
To handle this one needs the canonical
map
\begin{equation}\label{pi}
		  \pi : \mathrm{End}_k(A) \rightarrow A,
		  \quad \varphi \mapsto \varphi(1),
\end{equation} 
and the fact that 
$\mathrm{End}_k(A)$
is an algebra over 
$A^e$, with 
$$
		  \eta_{\mathrm{End}_k(A)} :
		  A^e \rightarrow \mathrm{End}_k(A),\quad
		  (\eta_{\mathrm{End}_k(A)}
		  (a \otimes b))(c):=acb.
$$
Now it makes sense to require 
$ \Delta $ and 
$ \varepsilon $ to factor through 
$ \iota $ and $ \pi $:
\begin{dfn}\label{takeke}
A (left) $\times_A$-bialgebra is an 
algebra $U$ over $A^e$ together with
two homomorphisms 
$\hat \Delta : U \rightarrow U \times_A U$
and $\hat \varepsilon : U \rightarrow 
\mathrm{End}_k(A)$
of algebras over $A^e$ 
such that 
$U$ is a coalgebra in $\amoda$ via
$\Delta = \iota \circ \hat \Delta$
and $ \varepsilon = \pi \circ \hat
 \varepsilon$. 
\end{dfn}

So one has for example for any $\times_A$-bialgebra 
$$
		  \Delta (a \lact u \ract b)=
		  a \lact u_{(1)} \otimes_A 
		  u_{(2)} \ract b,\quad
		  \Delta (a \blact u \bract b)=
		  u_{(1)} \bract b \otimes_A a \blact u_{(2)},
$$
where we started to use Sweedler's shorthand 
notation 
$u_{(1)} \otimes_A u_{(2)}$ 
for $\Delta(u)$.

Be aware that the four
$A$-actions are not the only feature
of $\times_A$-bialgebras that disappears
for $A=k$. Another 
crucial one is for example that
the counit $ \varepsilon : U \rightarrow A$ is 
not necessarily a ring homomorphism.
Note also that many authors write 
$s(a):=\eta (a \otimes 1)$ and
$t(a):=\eta (1 \otimes a)$ and formulate
the theory using these so-called source and 
target maps rather than $ \eta $.

\subsection{The monoidal category
  $\umod$ \cite{Schau:BONCRAASTFHB}} 
Definition~\ref{takeke}
might appear complicated, but it is 
the correct concept from several points 
of view. An
important one for us is the following result
of Schauenburg \cite[Theorem~5.1]{Schau:BONCRAASTFHB}:
\begin{thm}\label{schaui1}
The $\times_A$-bialgebra 
structures on an algebra 
$\eta : A^e \rightarrow U$ over 
$A^e$ correspond
bijectively to monoidal structures 
on $\umod$ for which
the forgetful functor 
$\umod \rightarrow \amoda$ induced by 
$\eta$ is strictly monoidal.  
\end{thm}

Given a $\times_A$-bialgebra structure
on $U$, the monoidal structure on
$\umod$ is defined as for bialgebras: 
one takes the tensor product 
$M \otimes_A N$ of 
the $A$-bimodules underlying 
$M,N \in \umod$ 
and defines a left
$U$-action via $\Delta$, 
\begin{equation}\label{roaig}
		  u(m \otimes_A n):=u_{(1)}m
		  \otimes_A
		  u_{(2)}n,\quad u \in U,m \in M,n
		  \in N.
\end{equation} 
\begin{dfn}
If $U$ is a $\times_A$-bialgebra and 
$M,N \in \umod$ are left $U$-modules, we
denote the left $U$-module 
$M \otimes_A N$ with
$U$-action (\ref{roaig}) by 
$M \otimes N$. 
\end{dfn} 
The unit object in $\umod$ is $A$ 
on which $U$ acts
via
$$
		  \hat \varepsilon (u)(a)=
		  \varepsilon(a \blact u)=
		  \varepsilon (u \bract a),
$$ 
where the last equality  
is a consequence of the definition 
of a $\times_A$-bialgebra.

There is an analogous notion of 
right $\times_A$-bialgebra for which 
$\modu$ is monoidal. 
However, 
for a left $\times_A$-bialgebra there is
in general  
no canonical monoidal structure on 
$\modu$ or even only 
right action of $U$ on $A$.
 
\subsection{$\times_A$-Hopf algebras \cite{Schau:DADOQGHA}} 
Let $U$ be a
$\times_A$-bialgebra and define
\begin{equation}
\label{Galois}
\beta: {}_\blact U \otimes_\Aop U_\ract \to 
U_\ract \otimes_A {}_\lact U, \quad u
\otimes_\Aop v \mapsto  
u_{(1)}  \otimes_A u_{(2)}v, 
\end{equation}
the so-called Galois map of $U$, where 
$$
 {}_\blact U \otimes_\Aop U_\ract = U
 \otimes_k U/
{{\rm
 span}\{a \blact u \otimes_k v - u
 \otimes_k v \ract a\,|\,
u,v \in U, a \in A \}}.
$$
One could flip the tensor components
in order to avoid taking the tensor
product over $A^\mathrm{op}$, but we
found it more convenient to
keep $\beta$ in the form which is
standard for bialgebras over fields. For
the latter it is
 easily seen that $\beta$ is bijective if
and only if $U$ is a Hopf algebra with
$\beta^{-1}(u \otimes_k v) := u_{(1)}
\otimes S(u_{(2)})v$, where $S$ is the
antipode of $U$. This motivates
the following definition due to 
Schauenburg~\cite{Schau:DADOQGHA}:
\begin{dfn}
A $\times_A$-bialgebra $U$ is 
a $\times_A$-Hopf  algebra if
$\beta$ is a bijection.
\end{dfn}

Following Schauenburg, we adopt a Sweedler-type notation  
\begin{equation}
\label{pm}
u_+ \otimes_\Aop u_- := \beta^{-1}( u \otimes_A 1)
\end{equation}
for the so-called translation map
$$
		  \beta^{-1}(\cdot \otimes_A 1) : U \rightarrow 
		  {}_\blact U \otimes_\Aop U_\ract.
$$ 
Since substantial for the subsequent
calculations, we 
list some 
properties of $\beta^{-1}$ as proven in
\cite[Proposition~3.7]{Schau:DADOQGHA}: one 
has for all $u, v \in U$, $a \in A$
\begin{eqnarray}
\label{Sch1}
u_{+(1)} \otimes_A u_{+(2)} u_- &=& u \otimes_A 1 \in U_\ract \otimes_A {}_\lact U \\
\label{Sch2}
u_{(1)+} \otimes_\Aop u_{(1)-} u_{(2)}  &=& u \otimes_\Aop 1 \in  {}_\blact U
\otimes_\Aop U_\ract \\ 
\label{Sch3}
u_+ \otimes_\Aop u_- & \in 
& U \times_\Aop U\\
\label{Sch37}
u_+ \otimes_\Aop u_{-(1)} \otimes_A
u_{-(2)} 
&=& u_{++}
\otimes_\Aop u_-
\otimes_A u_{+-}\\
\label{Sch4}
(uv)_+ \otimes_\Aop (uv)_- &=& u_+v_+
\otimes_\Aop v_-u_- 
\\ 
\label{Sch5}
\eta(a \otimes b)_+ \otimes_\Aop 
\eta(a \otimes b)_- &=& 
\eta(a \otimes 1) \otimes_\Aop 
\eta(b \otimes 1),
\end{eqnarray}
where in (\ref{Sch3}) we abbreviated
$$
		  U \times_\Aop U:=
		  \Bigl\{\sum_i u_i \otimes_\Aop v_i \in 
		  {}_\blact U \otimes_\Aop U_\ract\,|\,
		  \sum_i u_i \ract a \otimes_\Aop v_i=
		  \sum_i u_i \otimes_\Aop a \blact
		  v_i
		  \Bigr\}
$$
and in (\ref{Sch37}) the tensor product
over $A^\mathrm{op}$ links the first and
third tensor component (cf.~\cite
[Equation (3.7)]{Schau:DADOQGHA}).
By (\ref{Sch1}) and (\ref{Sch3})
one can write 
\begin{equation}
\label{pmb}
\beta^{-1}(u \otimes_A v) = u_+ \otimes_\Aop u_-v
\end{equation}
 which is easily checked to be
 well-defined over $A$ with (\ref{Sch4}) 
 and (\ref{Sch5}).

\subsection{Examples}\label{secex} 
Clearly, Hopf algebras over $k$ 
such as universal enveloping algebras of
Lie algebras or group algebras
are $\times_k$-Hopf
algebras. But also the
enveloping algebra of an associative
algebra that governs Hochschild
(co)homology is an example as 
pointed out by 
Schauenburg \cite{Schau:DADOQGHA}: 
\begin{ex}
The enveloping
 algebra $U:=A^e$
of any $k$-algebra $A$ is a
$\times_A$-bialgebra 
with $\eta=\mathrm{id}_{A^e}$ and 
coproduct and counit
$$
\Delta: U \to U \otimes U, \> a
 \otimes_k b \mapsto (a \otimes_k 1) 
\otimes_A 
(1 \otimes_k b),\quad 
\varepsilon: U \to A, \> a \otimes_k
 b \mapsto ab.
$$
As for the $\times_A$-Hopf algebra
 structure, the tensor product in 
question reads
$$
 {}_\blact U \otimes_\Aop U_\ract = U
 \otimes_k U/{\rm span}_k\{(a \otimes_k cb) 
\otimes_k (a' \otimes_k b')
 -(a \otimes_k b) \otimes_k 
(a' \otimes_k b'c)\},
$$
where $cb$ and $b'c$ is understood to be
the product in $A$.
One then easily verifies that
$$
(a \otimes_k b)_+ \otimes_\Aop 
(a \otimes_k b)_- := (a \otimes_k 1)
\otimes_\Aop (b \otimes_k 1)
$$
yields an inverse of the Galois map 
defined as in (\ref{pmb}).
\end{ex}

Finally we discuss
Lie-Rinehart algebras which 
define for example
Poisson (co)homology. Several authors 
\cite{Xu:QG,KhaRan:CCOEHA,MoeMrc:OTUEAOALRA} have shown that
their enveloping algebras are 
$\times_A$-bialgebras,
but they are 
in fact also $\times_A$-Hopf algebras:  

\begin{ex}\label{bumpy}
Let $(A,L)$ be a Lie-Rinehart algebra
over $k$ \cite{Rin:DFOGCA,Hue:PCAQ}. We 
denote by $(a,X) \mapsto aX$ the
$A$-module structure on $L$ and
by $(X,a) \mapsto X(a)$ the 
 $L$-action on $A$ given by the anchor
 $\hat \varepsilon : L \to \mathrm{Der}_k(A)$. 
Its universal enveloping
algebra $U=U(A,L)$ is the
universal $k$-algebra equipped with
two maps 
$$
		  \iota_A: A \to U,\quad \iota_L: L \to U
$$
of $k$-algebras and of $k$-Lie algebras,
respectively, and subject to the identities
$$
		  \iota_A(a) \iota_L(X) = \iota_L(aX),\quad 
		  \iota_L(X) \iota_A(a) - \iota_A(a) \iota_L(X) =
		  \iota_A(X(a))
$$
for $a \in A,X \in L$; confer
\cite{Rin:DFOGCA} for the precise 
construction.
The map
 $\iota_A$ is injective, so we refrain 
from further mentioning it. We 
will also merely write $X$ when we 
mean $\iota_L(X)$ (if $L$ is
 $A$-projective, then $\iota_L$ is injective
 as well).

Recall now from e.g.\
 \cite{Xu:QG,MoeMrc:OTUEAOALRA}  
that $U$ carries the
structure of a
$\times_A$-bialgebra: the maps
 $\eta(- \otimes 1)$ and $\eta(1 \otimes
 -)$ 
are equal and given by $\iota_A$. The
prescriptions 
\begin{equation}
\label{LCoProdVL1}
\Delta (X) = 1 \otimes_A X + X \otimes_A
1, \quad \Delta (a)
= a \otimes_A 1
\end{equation}
which map $X \in L$ and 
$a \in A$ into $U \times_A U$ can
be extended by the universal
property 
to a coproduct $\hat\Delta: U \to U \times_A U$. 
The counit is similarly 
given by the extension of the 
anchor $\hat \varepsilon$ to $U$.
The bijectivity of the Galois
map is seen in the same way: 
the translation map is given on
generators as
\begin{equation}
\label{LieRineIsHopfX}
a_+ \otimes_{\Aop} a_- := 
a \otimes_\Aop 1,\quad 
X_+ \otimes_\Aop X_- :=
X \otimes_\Aop 1 - 1 \otimes_\Aop X.
\end{equation} 
These maps stay in 
 $U \times_\Aop U$ which is an algebra 
through the product of $U$ in the
 first and its opposite  
in the second tensor 
 factor. By universality
 we obtain a map $U \to U \times_\Aop U$,
 and then $\beta^{-1}$ is defined using
 (\ref{pmb}).

On the other hand, $U$ is not
necessarily a Hopf algebroid in the
sense of \cite{Boe:HA} (this also
 answers B\"ohm's question asked therein
 whether any $\times_A$-Hopf algebra is a Hopf
 algebroid). This structure
assumes the existence of an antipode 
$S : U \rightarrow U^\mathrm{op}$
satisfying certain axioms. As a result, 
the left $U$-action on $A$ yields by
composition with $S$ also a right
 $U$-module structure. However, there
might be an obstruction for this. For
example, take
$L=\Gamma(T^{1,0}S^2)$, where 
$T^{1,0}S^2 \oplus T^{0,1}S^2=TS^2
 \otimes \mathbb{C}$ is the
 decomposition of the complexified
 tangent bundle of the 2-sphere $S^2
 \subset \mathbb{R}^3$ into the
 holomorphic and antiholomorphic part
 with respect to the standard complex
 structure. This defines together with 
$A=C^\infty(S^2,\mathbb{C})$ a
Lie-Rinehart algebra, where the action
of $L$ on $A$ is the usual action of a
vector field on a smooth function and
the action of $A$ on $L$ is given by
fibrewise multiplication. We know by work
of Huebschmann \cite{Hue:LRAGAABVA}
that the right $U$-module structures on
$A$ correspond bijectively to left
$U$-module structures on 
$L$ itself (in general on its top
exterior power  over $A$, but here this
is $L$ because $T^{1,0}S^2$ is only a
line bundle). 
Such a left $U$-action
corresponds precisely to a flat
connection $\nabla$ on the complex line bundle 
$T^{1,0}S^2$, with $X \in L$ acting on
 sections of $T^{1,0}S^2$ by the covariant
 derivative $\nabla_X$ (see 
\cite{Hue:LRAGAABVA} for the details). 
But the 
curvature of any connection
represents the first Chern
class of the bundle which is nonvanishing
since $T^{1,0}S^2$ is not
trivial. Therefore, there is no flat
connection aka left $U$-action on $L$ and hence no right 
$U$-action on $A$.
\end{ex}

\section{Multiplicative structures}
\subsection{$\D^-(U)$ as a suspended 
monoidal category \cite{Sua:THHAFTACOCP} }
For any ring $U$, we denote by 
$\D^{-}(U)$ the derived category of 
bounded above cochain
complexes of left $U$-modules. As
usual, we identify any 
$M \in \umod$ with a complex in
$\D^-(U)$ concentrated in degree $0$, and
any bounded below chain complex 
$P_\bullet$ with a
bounded above cochain complex by putting 
$P^n:=P_{-n}$.

If $U$ is an 
$A$-biprojective $\times_A$-bialgebra, 
then any projective 
$P \in \umod$ is $A$-biprojective. Hence 
the monoidal structure of $\umod$
extends to a monoidal 
structure on $\D^-(U)$ 
with unit
object still given by $A$ and product 
being the total tensor product 
$\lotimes=\lotimes_A$ (the
$A$-biprojectivity of $U$-projectives is
needed for example to have 
\cite[Lemma~10.6.2]{Wei:AITHA}).

Together with the shift functor 
$T : \D^-(U) \rightarrow \D^-(U)$,
$(TC)^n=C^{n+1}$, $\D^-(U)$ becomes what
is called a suspended monoidal category 
in \cite{Sua:THHAFTACOCP}. This just means 
that for all $C,D \in \D^-(U)$, the canonical
isomorphisms 
$$
		  TC \lotimes D \simeq T(C \lotimes D)
		  \simeq C \lotimes TD
$$ 
given by the obvious renumbering make  
the diagrams
$$
\xymatrix{A \lotimes TC \ar[d]\ar[r] & TC
\ar[dl]\\
T(A \lotimes C)}\quad\quad
\xymatrix{TC \lotimes A \ar[d]\ar[r] & TC
\ar[dl]\\
T(C \lotimes A)}
$$ 
commutative and the diagram
$$
\xymatrix{TC \lotimes TD \ar[d]\ar[r] &
T(C \lotimes TD)
\ar[d]\\
T(TC \lotimes D) \ar[r] & T^2(C \lotimes D)}
$$
anticommutative (commutative up to a sign
$-1$). 

\subsection{$\smallsmile$ and 
$\circ$ \cite{Sua:THHAFTACOCP} }\label{absch1}   
As a special case of the constructions
from \cite{Sua:THHAFTACOCP}, we 
define 
for any $A$-biprojective 
$\times_A$-bialgebra $U$
and $L,M,N \in \umod$ the cup product
$$
		  \smallsmile \,\,: 
		  \mathrm{Ext}^m_U(A,M) \times 
		  \mathrm{Ext}^n_U(A,N)
		  \rightarrow 
		  \mathrm{Ext}^{m+n}_U(A,M \otimes
		  N)
$$ 
and the classical Yoneda product
$$
		  \circ : 
		  \mathrm{Ext}^m_U(N,M) \times 
		  \mathrm{Ext}^n_U(L,N)
		  \rightarrow 
		  \mathrm{Ext}^{m+n}_U(L,M).
$$ 
The latter is just the composition 
of morphisms in $\D^-(U)$ if one
identifies 
$$
		  \mathrm{Ext}_U^n(L,N)
		  \simeq \mathrm{Hom}_{\D^-(U)}(L,T^nN),
$$
and
$$
		  \mathrm{Ext}_U^m(N,M)
		  \simeq \mathrm{Hom}_{\D^-(U)}(N,T^mM)
		  \simeq
		  \mathrm{Hom}_{\D^-(U)}(T^nN,T^{m+n}M).
$$
The former is obtained as follows: given
$$
		  \varphi \in \mathrm{Ext}_U^m(A,M)
		  \simeq
		  \mathrm{Hom}_{\D^-(U)}(A,T^mM),
$$
$$
		  \psi \in \mathrm{Ext}_U^n(A,N)
		  \simeq \mathrm{Hom}_{\D^-(U)}(A,T^nN),
$$
one defines $\varphi \smallsmile \psi$
as the composition
\begin{eqnarray}
&& A \simeq  
		 A \otimes A \nonumber\\ 
&\xymatrix{\ar[r]^{\varphi
		 \otimes \psi} &}& 
 		 T^mM \lotimes T^nN \simeq T^m(M
		 \lotimes T^nN) \simeq
		 T^{m+n}(M \lotimes N)
		 \nonumber\\ 
&\xymatrix{\ar[r] &}& T^{m+n}(M \otimes N),
		 \nonumber
\end{eqnarray} 
where the last map is the augmentation 
$M \lotimes N \rightarrow H^0(M \lotimes
N) \simeq \mathrm{Tor}^A_0(M,N) \simeq M
\otimes N$, or rather $T^{m+n}$ applied
to this morphism in $\D^-(U)$.

A straightforward extension of
Theorem~{1.7} from 
\cite{Sua:THHAFTACOCP} now gives:
\begin{thm}\label{suarez}
If $U$ is an $A$-biprojective 
$\times_A$-bialgebra,
then we have
$$
		  \psi \circ \varphi = \varphi
		  \smallsmile \psi = (-1)^{mn}\psi \smallsmile
		  \varphi,\quad
		  \varphi \in
		  \mathrm{Ext}_U^m(A,A),
		  \psi \in \mathrm{Ext}_U^n(A,M)
$$
as elements of 
$
		  \mathrm{Ext}^{m+n}_U(A,M) \simeq
		  \mathrm{Ext}_U^{m+n}(A,A \otimes
		  M) \simeq
 \mathrm{Ext}_U^{m+n}(A,M \otimes A).
$

In particular, $\mathrm{Ext}_U(A,A)$
 becomes through either of the products
 a graded commutative algebra over 
the commutative subring $\mathrm{Hom}_U(A,A)$.
\end{thm}
\begin{pf}
This is proven exactly as in
 \cite{Sua:THHAFTACOCP}. 
For the reader's convenience we include
one of the diagrams involved. The
 unlabeled arrows are 
canonical maps coming from
the suspended monoidal structure.
$$
\xymatrix@=38pt{
A \ar[d]_\varphi \ar[r]
& A \otimes A \ar[d]_{\mathrm{id}
 \otimes \varphi} \ar[rd]^{\psi \otimes
 \mathrm{id}}\\
T^m A \ar[dr] \ar[r]\ar[dd]_{\mathrm{id}} & A \otimes
 T^m A \ar[d] \ar[dr]^{\psi \otimes
 \mathrm{id}}& T^n M
 \lotimes A \ar[d]^{\mathrm{id} \otimes
 \varphi}\\
& T^m(A \otimes A) \ar[dl]
 \ar[d]_{T^m(\psi \otimes \mathrm{id})\!\!}& 
T^n M \lotimes T^m A \ar[d] \ar[dl]\\
T^m A \ar[d]_{T^m(\psi)} & T^m(T^n M
 \lotimes A) \ar[dl] \ar[d]& T^n(M
 \lotimes T^mA)\ar[dl]\\
T^{m+n}M & T^{m+n}(M \lotimes A) \ar[l]
}
$$
The morphism 
$\psi \circ \varphi \in \mathrm{Hom}_{\D^-(U)}(A,T^{m+n}M)$
is the path going straight down 
from $A$ to $T^{m+n}M$, and $\psi
\smallsmile \varphi$ is the one which
goes clockwise round the whole diagram. All
faces of the 
diagram commute except the
 lower right square
which introduces a sign $(-1)^{mn}$,
so we get  
$\psi \circ \varphi = (-1)^{mn} \psi
 \smallsmile \varphi$. The other
 identity is shown with a similar diagram.
\end{pf}
{}\\[-10mm]

\subsection{Tensoring 
projectives}\label{barreso} 
This paragraph is a small excursus about
the projectivity of the tensor product of two 
projective objects of a monoidal
category.
For example, $U \otimes U \in \umod$ 
is not necessarily projective even for a
bialgebra $U$ over a field $A=k$
(so the $A$-projectivity of $U$ or the
exactness of $\otimes$ does not
help). Here is a simple example 
(for a detailed study of examples 
of categories
of Mackey functors see 
\cite{Gau:WPDNIFAOHA}): 
\begin{ex}\label{bsp}
Consider the bialgebra 
$U=\mathbb{C}[a,b,c]$ over
$A=k=\mathbb{C}$ with 
$$
		  \Delta (a)=a \otimes a,\quad 
		  \Delta (b)=a \otimes b+b \otimes c,\quad
		  \Delta (c) = c \otimes c,
$$
$$
		  \varepsilon (a) =1,\quad 
		  \varepsilon (b) = 0,\quad
		  \varepsilon (c) = 1.
$$
Geometrically, this is the coordinate ring 
of the complex algebraic semigroup 
$G$ of upper triangular $2 \times 2$-matrices, 
and $ \Delta $ and $ \varepsilon $ are dual 
to the semigroup law $G \times G \rightarrow G$
and the embedding of the
identity matrix into $G$.

We prove that 
$U \otimes U \in \umod$ 
is not projective by 
considering the fibres
of the semigroup law $G \times G \rightarrow G$. 
The fibre over a generic and hence invertible
element is 3-dimensional, but over $0$ it is 
 4-dimensional, and this will imply our
claim. We can use for example \cite[Theorem 19 on p.~79]{Mat:CA}:
\begin{thm}
Let $U \subset V$ 
be a flat extension of
commutative Noetherian rings, 
$\mathfrak{p} \subset V$ be a prime
 ideal and
$\mathfrak{q} := 
U \cap \mathfrak{p}$. Then 
$$
		  \mathrm{dim} (V_\mathfrak{p}) = 
		  \mathrm{dim} (U_\mathfrak{q}) + 
		  \mathrm{dim} (V_\mathfrak{p} \otimes_U
		  U(\mathfrak{q})),
$$ 
where $ \mathrm{dim} $ denotes the Krull 
dimension of a ring, $V_\mathfrak{p}$ is 
the localisation of $V$ at
$\mathfrak{p}$ and  
$U(\mathfrak{q}):=U_\mathfrak{q} / 
\mathfrak{q} U_\mathfrak{q}$ 
is the 
residue field of the localisation
$U_\mathfrak{q}$. 
\end{thm}

Apply this to our example 
$U \simeq \Delta (U) \subset V:=U
 \otimes U$: 
let $\mathfrak{p}$ be the ideal of $V$ 
generated by $a \otimes_\mathbb{C} 1$,
$1 \otimes_\mathbb{C} a$, $b \otimes_\mathbb{C} 1$,
$1 \otimes_\mathbb{C} b$, $c \otimes_\mathbb{C} 1$, 
$1 \otimes_\mathbb{C} c$. Geometrically, $V$ is 
the coordinate ring of $\mathbb{C}^6$ and 
$V_\mathfrak{p}$ is the local ring in 
$0$, so $\mathrm{dim} (V_\mathfrak{p})=6$.
Since $1 \notin \mathfrak{p}$, 
$\mathfrak{q}=U \cap \mathfrak{p}$
is proper, and it contains the ideal 
generated by 
$\Delta(a)=a \otimes_\mathbb{C} a$,
$\Delta(b)=a \otimes_\mathbb{C} b+b \otimes_\mathbb{C} c$,
$\Delta (c)=c \otimes_\mathbb{C} c$ which is maximal in 
$U$, so
$\mathfrak{q} \subset U$ is the ideal generated
by $a,b,c$, and $U_\mathfrak{q}$
is the local ring of 
$\mathbb{C}^3$ at $0$ with 
$\mathrm{dim} (U_\mathfrak{q})=3$.
The field $U(\mathfrak{q})$ is 
obviously $\mathbb{C}$, and we can write 
$V_\mathfrak{p} \otimes_U U(\mathfrak{q})$
also as $V_\mathfrak{p}/
\Delta(\mathfrak{q})V_\mathfrak{p}$.
Since $ \Delta (\mathfrak{q}) V_\mathfrak{p}$
is contained in the ideal $\mathfrak{r}$
generated in $V_\mathfrak{p}$ by the 
elements $a \otimes_\mathbb{C} 1,1 \otimes_\mathbb{C} c$,
we have 
$\mathrm{dim} (V_\mathfrak{p}/
\Delta(\mathfrak{q})V_\mathfrak{p})
\ge
\mathrm{dim} (V_\mathfrak{p} /\mathfrak{r})$. 
Now $V_\mathfrak{p}/\mathfrak{r}$
is the local ring of 
$\mathbb{C}^4 \subset \mathbb{C}^6$
at $0$ and hence
$\mathrm{dim} 
(V_\mathfrak{p} /\mathfrak{r})=4$.
In total, we get the strict inequality 
$3+\mathrm{dim} (V_\mathfrak{p}/
\Delta(\mathfrak{q})V_\mathfrak{p})
\ge 3+4=7>6$, 
and hence $V$ is not flat over 
$U$ and in particular not projective.   
\end{ex}

For $\times_A$-Hopf algebras the
situation is, however, much simpler: 
notice that
$$
{}_\blact U \otimes_\Aop
M \ract := U \otimes_k M/{{\rm
 span}\{a \blact u \otimes_k m - 
u \otimes_k m \ract a\,|\,u \in U,
a \in A,m \in M\}}
$$ 
is for any $\times_A$-bialgebra $U$ and 
$M \in \umod$ a left $U$-module by
left multiplication 
on the first factor. Just as for $M=U$, 
there is a Galois map 
$$
		  \beta_M : {}_\blact U \otimes_\Aop 
		  M_\ract \to  U \otimes M, \quad 
		  u \otimes_\Aop m \mapsto u_{(1)} \otimes_A
     	  u_{(2)} m,
$$
and we have:
\begin{lem}
\label{resomod}
For any $\times_A$-bialgebra $U$, the 
generalised Galois map $\beta_M$ is a
morphism of $U$-modules. If $U$ is a
 $\times_A$-Hopf algebra, then 
$\beta_M$ is bijective.
\end{lem}
\begin{pf}
The $U$-linearity of $\beta_M$ follows
immediately from the fact that 
$\hat\Delta : U \rightarrow U \times_A U$
is a homomorphism of algebras over $A^e$.
Furthermore, if $\beta$ is a bijection,
 then $\beta_M$ is so as well 
since we can identify $\beta_M$ with 
$\beta \otimes_U {\rm id}_M$, and then
the inverse is simply
given by
$\beta_M^{-1}(u \otimes_A m) = u_+
 \otimes_\Aop u_-m$. 
\end{pf}
{}\\[-7mm]

Using this one now gets:
\begin{thm}\label{uuu}
If $U$ is a $\times_A$-Hopf algebra
and $U_\ract \in \moda$ is projective, then 
$P \otimes Q \in \umod$ is
projective for all projectives $P,Q \in
 \umod$. 
\end{thm} 
\begin{pf}
By assumption, any projective module 
over $U$ is also projective over
 $A^\mathrm{op}$, and
if $\varphi : R \rightarrow S$ is 
any ring map, then 
$S \otimes_R \cdot : \rmod \rightarrow \smod$
maps projectives to projectives. 
This shows that
${}_\blact U \otimes_\Aop U_\ract$ 
and hence (Lemma~\ref{resomod})
$U \otimes U$ is projective. 
Since $\otimes=\otimes_A$ commutes
with arbitrary direct sums, 
$P \otimes Q$ is projective for 
all projectives $P,Q$. 
\end{pf}
{}\\[-9mm]
\begin{cor}
If $U$ is as in Theorem~\ref{uuu} and
$P \in \D^-(U)$ is a
projective resolution of $A \in \umod$, 
then so is $P \otimes P:=\mathrm{Tot}(P_\bullet \otimes
 P_\bullet)=P \lotimes P$.
\end{cor}

This leads to the traditional 
construction of $\smallsmile\,$
given for $A=k$ in 
\cite[Chapter XI]{CarEil:HA}: 
one fixes a projective resolution
$P$ of $A$, and by the above,  
$ \mathrm{Ext}_U(A,M \otimes N)$
is the total (co)homology of the double (cochain)
complex
$$
C^2_{mn}:=
\mathrm{Hom}_U(P_m \otimes P_n,
M \otimes N).
$$
Then $\smallsmile$ is given as the composition of the canonical map 
\begin{eqnarray}
&& \bigoplus_{m+n=p} 
		  \mathrm{Ext}^m_U(A,M) \otimes_k
		  \mathrm{Ext}^n_U(A,N)\nonumber\\  
&\simeq&
		  \bigoplus_{m+n=p} 
		  H^m(\mathrm{Hom}_A(P_\bullet,M))
		  \otimes_k 		  
 H^n(\mathrm{Hom}_A(P_\bullet,N))\nonumber\\ 
&\rightarrow& H^p(\bigoplus_{m+n=\bullet} 
		  \mathrm{Hom}_A(P_m,M)
		  \otimes_k 		  
 \mathrm{Hom}_A(P_n,N))
= H^p(\mathrm{Tot}(C^1_{\bullet\bullet}))
		  \nonumber 
\end{eqnarray} 
where $C^1_{mn}:=\mathrm{Hom}_U(P_m,M) 
\otimes_k \mathrm{Hom}_U(P_n,N)$, with the map  
$$
		  H(\mathrm{Tot}(C^1_{\bullet\bullet})) \rightarrow 
		  H(\mathrm{Tot}(C^2_{\bullet\bullet}))
		  \simeq
		  \mathrm{Ext}_U(A,M \otimes N)
$$
that is induced by 
the morphism of double complexes
$$
		  C^1_{mn} \ni \varphi \otimes_k 
		  \psi \mapsto 
		  \{x \otimes y \mapsto 
		  \varphi(x) \otimes \psi(y)\}
		  \in C^2_{mn}.
$$

For the sake of completeness let us
finally remark that as for $A=k$ one can in particular use 
the bar construction 
to obtain a canonical resolution:
\begin{lemma}
For any $\times_A$-bialgebra $U$,
the complex of left $U$-modules 
$$
		  \Cbar_n:=
		  ({}_\blact
		  U_\ract)^{\otimes_\Aop n+1},\quad
		  u(v_0 \otimes_\Aop \cdots \otimes_\Aop
		 v_n):=uv_0 \otimes_\Aop \cdots \otimes_\Aop v_n 
$$
whose boundary map is given by
\begin{eqnarray}
b' : u_0 \otimes_\Aop \cdots \otimes_\Aop u_n 
&\mapsto& \sum_{i=0}^{n-1} (-1)^i
 		  u_0 \otimes_\Aop \cdots \otimes_\Aop
		  u_iu_{i+1} \otimes_\Aop \cdots \otimes_\Aop
		  u_n \nonumber\\ 
&& + (-1)^n u_0 \otimes_\Aop \cdots 
		  \otimes_\Aop \varepsilon(u_n)
		  \blact u_{n-1}
		  \nonumber
\end{eqnarray} 
is a contractible resolution of $A \in \umod$ with
augmentation 
$$
		  \varepsilon : \Cbar_0=U \rightarrow 
		  A=:\Cbar_{-1},
$$ 
and if $U_\ract \in \moda$ is
 projective, then $\Cbar_n \in \umod$
is projective.
\end{lemma}
\begin{pf}
All claims are straightforward: 
there is a contracting homotopy
$$
		  s : \Cbar_n \rightarrow \Cbar_{n+1},\quad
		  u_0 \otimes_\Aop \cdots \otimes_\Aop u_n 
		  \mapsto
		  1 \otimes_\Aop u_0 \otimes_\Aop \cdots 
		  \otimes_\Aop u_n,\quad n \ge 0,
$$ 
$$
		  s : A=\Cbar_{-1} \rightarrow U=\Cbar_0,\quad
		  a
		  \mapsto
		  \eta (a \otimes 1),
$$
and the
projectivity of $\Cbar_n$
follows as in the proof of 
Theorem~\ref{uuu}. 
\end{pf}
{}\\[-10mm]

\subsection{The functor 
$\otimes : \umod \times \modu
  \rightarrow \modu$} 
Now we introduce the functor $\otimes$
mentioned in Theorem~\ref{main}.
\begin{lem}
\label{RechtsLinks}
Let $U$ be a $\times_A$-Hopf algebra and
$M \in \umod$, $P \in \modu$ be left and right $U$-modules, respectively. Then the formula
\begin{equation}
\label{pfeffer}
(m \otimes_A p)u:= u_-m \otimes_A pu_+, \quad\quad u \in U, \ m
\in M, p \in P,
\end{equation}
defines a right $U$-module structure on
the tensor product
\begin{equation}
\label{TensLR}
M \otimes_A P := M \otimes_k P/{{\rm
 span}\{m \ract a \otimes_k p - 
m \otimes_k a \blact p\,|\,
a \in A\}}.
\end{equation}
If $N$ is any other (left) $U$-module,
then the canonical isomorphism
\begin{equation}\label{luton}
		  (M \otimes N) \otimes_A P
			\simeq
			M \otimes_A (N \otimes_A P)
\end{equation} 
of $A$-bimodules is also an isomorphism 
in $\modu$. Finally, the tensor flip
$$
		  (M \otimes_A P) \otimes_U N
		  \rightarrow 
		  P \otimes_U (N \otimes_A M),\quad
		  m \otimes_A p \otimes_U n
		  \mapsto
		  p \otimes_U n \otimes_A m
$$ 
is an isomorphism of $k$-modules.
\end{lem} 
\begin{pf}
To show firstly that 
(\ref{pfeffer}) is well-defined
over $A$, we compute
\begin{equation*}
\begin{split}
\big(m \otimes_A (a \blact p)\big)u &= 
u_-m \otimes_A p\eta(1 \otimes a)u_+ 
= u_-m \otimes_A p(u_+ \ract a) \\
&= (a \blact u_-)m  \otimes_A pu_+
= u_- \big(\eta(1 \otimes a)m\big)
 \otimes_A pu_+  \\
&= \big((m \ract a) \otimes_A p\big)u,
\end{split}
\end{equation*}
where (\ref{Sch3}) and the action
 properties were used.  
Together with (\ref{TensLR}) this also proves the
well-definedness with respect to the presentation 
of $ u_+ \otimes_\Aop u_-$.
With the help of (\ref{Sch4}) one sees
immediately 
that for $u, v \in U$ we have 
$$
\big(m\otimes_A p\big)(uv) =  (uv)_-m \otimes_A p(uv)_+ =  v_-u_- m
\otimes_A pu_+v_+ 
= \big((m\otimes_A p)u\big)v,
$$
since $P$ and $M$ were right and left $U$-modules, respectively. 
As a conclusion, $M \otimes_A P \in \modu$.
Equation (\ref{luton}) is a direct
consequence of the associativity of the
tensor product of $A$-bimodules 
and of (\ref{Sch37}).

For the last part one has to check that the flip
is well-defined: we have
\begin{eqnarray}
 \eta (1 \otimes a)m \otimes_A p
		  \otimes_U n 
&\mapsto& p \otimes_U 
		  n \otimes_A \eta (1 \otimes a)m 
= p \otimes_U 
		  \eta (1 \otimes a)
		  (n \otimes_A m) \nonumber\\ 
&=& p \eta (1 \otimes a) \otimes_U
		  (n \otimes_A m), \nonumber
\end{eqnarray}
which is what $m \otimes_A p \eta (1 \otimes a)
 \otimes_U n$ 
gets mapped to. And secondly, we have
\begin{eqnarray}
 m \otimes_A p
		  \otimes_U u n 
&\mapsto& p \otimes_U 
		  u n \otimes_A m
= p \otimes_U (u_+)_{(1)} n \otimes_A
		  (u_+)_{(2)} u_- m \nonumber\\ 
&=& p \otimes_U u_+(n \otimes_A 
		  u_-m) 
= p u_+ \otimes_U n \otimes_A 
		  u_-m,\nonumber
\end{eqnarray} 
which is what 
$u_- m \otimes_A p u_+ \otimes_U
	  n = (m \otimes_A p)u \otimes_U
 n$ gets mapped to.
\end{pf}
{}\\[-9mm]
\begin{dfn}
We denote the above 
constructed $U^\mathrm{op}$-module 
by $M \otimes P$.
\end{dfn}

Thus an unadorned $\otimes$ refers
from now on either to the monoidal product on 
$\umod$ or to the just defined action of 
$\umod$ on $\modu$. For
example, (\ref{luton}) would now simply
be written as 
$ (M \otimes N) \otimes P \simeq
M \otimes (N \otimes P)$.
\begin{ex}
Let $(A,L)$ be a Lie-Rinehart algebra and $M$ be a left and $N$ a
right $U(A,L)$-module, respectively (or,
 in the terminology of \cite{Hue:PCAQ,
 Hue:DFLRAATMC}, left 
and right $(A,L)$-modules). Using
 (\ref{LieRineIsHopfX}), one gets the 
right $U(A,L)$-module
structure on $M \otimes_A N$
from formula (2.4) in \cite[p.~112]{Hue:DFLRAATMC}: 
$$
(m\otimes_A n)X = m \otimes_A nX - Xm
 \otimes_A n, \qquad m \in M, \ n \in N,
 \ X \in L. 
$$
\end{ex}

If we assume again
that $U$ is $A$-biprojective, then 
the above results extend directly to the derived
category $\D^-(U^\mathrm{op})$:
we obtain a functor 
$$
		  \lotimes = \lotimes_A : 
		  \D^-(U) \times
		  \D^-(U^\mathrm{op}) \rightarrow 
		  \D^-(U^\mathrm{op})
$$
and we have for all $M,N \in \D^-(U)$, $P \in
\D^-(U^\mathrm{op})$ canonical isomorphisms
\begin{equation}\label{weihnacht}
		  (M \lotimes N) \lotimes P \simeq 
		  M \lotimes (N \lotimes P),\quad
		  (M \lotimes P) \lotimes_U N
		  \simeq
		  P \lotimes_U (N \lotimes M).
\end{equation}

\subsection{$\smallfrown$ and $\bullet$}\label{absch2} 
These products are dual to $\smallsmile$
and $\circ$. The first one is
$$
		  \bullet : \mathrm{Ext}_U^m(L,M)
		  \times \mathrm{Tor}^U_n(N,L)
		  \rightarrow \mathrm{Tor}^U_{n-m}(N,M)
$$
which exists for
a ring $U$
and $L,M \in \umod$, $N \in \modu$:
an element
$$
		  \varphi \in
		  \mathrm{Ext}_U^m(L,M) \simeq
		  \mathrm{Hom}_{\D^-(U)}(L,T^mM)
$$
defines a morphism in $\D^-(\mathbb{Z})$
$$
		  N \lotimes_U L \rightarrow N
		  \lotimes_U T^mM,\quad
 		  x \otimes_U y \mapsto x \otimes_U
		  \varphi (y),
$$
and $\varphi \bullet \cdot$ is the
induced map in (co)homology
\begin{eqnarray}
&&\mathrm{Tor}_n^U(N,L) \simeq
   	H^{-n}(N \lotimes_U L) \nonumber\\
&\xymatrix{\ar[r]^{H^{-n}(\mathrm{id}
 \otimes \varphi)}&}& 
		  H^{-n}(N \lotimes_U T^mM) \simeq
		  H^{m-n}(N \lotimes_U M) \simeq
		  \mathrm{Tor}^U_{n-m}(N,M). \nonumber
\end{eqnarray}

For $M \in \umod,N \in \modu$ as before, the cap product 
$$
		  \smallfrown \,\,:
		  \mathrm{Ext}^m_U(A,M) \times
		  \mathrm{Tor}^U_n(N,A)
		  \rightarrow 
		  \mathrm{Tor}^U_{n-m}(M \otimes
		  N,A)		  
$$
involves the functor
$\otimes$ from the previous paragraph, 
so for this we want $U$ to be an $A$-biprojective $\times_A$-Hopf
algebra again. Similarly as for
$\bullet$, 
$$
		  \varphi \in
		  \mathrm{Ext}_U^m(A,M) \simeq
		  \mathrm{Hom}_{\D^-(U)}(A,T^mM)
$$
defines a morphism in $\D^-(k)$
\begin{eqnarray}
&& N \lotimes_U A \simeq
		  N \lotimes_U (A \otimes A)
		  \nonumber\\  
&\xymatrix{\ar[r]^{\mathrm{id} \otimes
 \mathrm{id} \otimes \varphi} &} & N \lotimes_U (A \lotimes T^mM) \simeq
		  N \lotimes_U (T^mA \lotimes M) \simeq
		  (M \lotimes N) \lotimes_U T^mA
		  \nonumber\\
&\xymatrix{\ar[r]&}& (M \otimes N)
 \lotimes_U T^mA,\nonumber
\end{eqnarray} 
where the last $\simeq$ in the second
line is induced
by the tensor flip as in the
derived version (\ref{weihnacht}) of
Lemma~\ref{RechtsLinks}, and the morphism from
the second to the third line is
similarly as in the definition of
$\smallsmile$ induced by the morphism 
$M \lotimes N \rightarrow M \otimes N$ 
in $\D^-(U^\mathrm{op})$ that takes
zeroth cohomology. Passing now to
cohomology we get
$\varphi \smallfrown \cdot :
\mathrm{Tor}^U_n(N,A) 
\rightarrow \mathrm{Tor}_{n-m} (M
\otimes N,A)$.

More explicitly, if $P \in \D^-(U)$ is a
projective resolution of $A$, 
then $\smallfrown$ is
induced by the morphism 
$$
		  B^1_{ij} \ni n \otimes_U 
		  (x \otimes_A y) \mapsto 
		  \{\varphi \mapsto 
		  (\varphi (y) \otimes_A n) \otimes_U
		  x\} \in B^2_{ij}
$$
from the double complex
$$
		  B^1_{ij}:=N \otimes_U 
		  (P_i \otimes_A P_j)
$$
whose total homology is 
$\mathrm{Tor}^U(N,A)$ to the
double complex 
$$
		  B^2_{ij}:=\mathrm{Hom}_k(
		  \mathrm{Hom}_U(P_j,M),
		  (M \otimes N) \otimes_U P_i)
$$
whose homology has a natural map to 
$\mathrm{Hom}_k(\mathrm{Ext}_U(A,M),
\mathrm{Tor}^U(M \otimes N,A))$.

In direct analogy with
Theorem~\ref{suarez} we get:
\begin{thm}\label{suarez2}
If $U$ is an $A$-biprojective 
$\times_A$-Hopf algebra,
then we have 
$$
		  \varphi \bullet (x \otimes_U y) = \varphi
		  \smallfrown (x \otimes_U y),\quad
		  \varphi \in
 \mathrm{Ext}_U^m(A,A),
x \otimes_U y \in N \lotimes_U A
$$
as elements of 
$
		  N \lotimes_U A \simeq
		  (A \otimes N) \lotimes_U A.
$
\end{thm}

\section{Duality and the proof of Theorem~\ref{main}} 
\subsection{The underived case} 
In the special case
that $A$ is finitely
generated projective itself, 
Theorem~\ref{main} reduces to 
standard linear algebra. We go through 
this case first since it is
both instructive and used in the
proof of the general case. 
For the reader's convenience
we include full proofs.

\begin{lemma}\label{cardiff}
Let $U$ be a ring, $A \in \umod$ be finitely
generated projective, and
$A^*$ be $\mathrm{Hom}_U(A,U)$ 
with its canonical
 $U^\mathrm{op}$-module structure.
\begin{enumerate}
\item[1.] $A^*$ is
finitely generated projective, and 
if 
$e_1,\ldots,e_n $ are generators of $A$, then there
exist generators
$e^1,\ldots,e^n \in A^*$
with
$$
		  \sum_i e^i(a)e_i=a,\quad
		  \sum_i e^i \alpha (e_i)=\alpha
$$ 
for all $a \in A$ and $\alpha \in A^*$.
The element
$$
		  \omega:=\sum_i e^i \otimes e_i \in A^*
		  \otimes_U A
$$
is independent of the choice of
the generators $e_i,e^j$. 
\item[2.] For all
$U^\mathrm{op}$-modules $M$,
the assignment 
$$
		  \delta (m \otimes a)(\alpha):=
		  m \alpha (a),\quad
		  m \in M,a \in A,\alpha \in
		  A^*
$$
extends uniquely to an 
isomorphism of abelian groups
$$
		  \delta : M
		  \otimes_U A \rightarrow 
		  \mathrm{Hom}_{U^\mathrm{op}}(A^*,M).
$$
\item[3.] One has $(A^*)^* \simeq A$ 
and $A^* \otimes_U M \simeq
\mathrm{Hom}_U(A,M)$
for $M \in \umod$.
\end{enumerate} 
\end{lemma}
\begin{pf}
Since $A$ is projective, there is a
splitting $\iota : A \rightarrow U^n$
of  
$$ \pi : U^n \rightarrow A,\quad 
(u_1,\ldots,u_n) \mapsto \sum_i
 u_ie_i.
$$ Hence 
$U^n  \simeq A \oplus A_\perp$ for some
$A_\perp \in \umod$. Dually this 
gives $A^* \oplus (A_\perp)^* =
(U^n)^* \simeq U^n$, whence 
$A^*$ is finitely generated projective. 
The $e^i$ can be defined as the
composition of $\iota$ with the
projection of $U^n$ on the $i$-th summand. 
This proves the first parts of 1. For
2.~just note that 
$$
		  \mathrm{Hom}_{U^\mathrm{op}}
		  (A^*,M) \ni 
		  \varphi \mapsto
		  \sum_i \varphi (e^i) \otimes e_i
		  \in M \otimes_U A
$$
inverts $\delta$. Since
$\omega =
\delta^{-1}(\mathrm{id}_{A^*})$,
it does indeed not
depend on the choice of generators.
3.~now follows from 1.~and 2.
\end{pf}
  
As in the introduction, let us  
abbreviate in the situation
of this theorem 
$$
		  H^0(M):=\mathrm{Hom}_U(A,M),\quad
		  H_0(N):=N \otimes_U
		  A
$$
for $M \in \umod$, 
$N \in \modu$, and call
$\omega \in H_0(A^*)$
the fundamental class
of $(U,A)$. Then 3.~says for $M=A$
that we have an isomorphism
\begin{equation}\label{sonn}
		  \cdot \bullet \omega : 
		  H^0(A) \rightarrow
		  H_0(A^*),\quad
		  \varphi \mapsto \sum_i
		  e^i \otimes \varphi (e_i).
\end{equation} 
Using Lemma~\ref{RechtsLinks} we can upgrade this
to the underived case of Theorem~\ref{main}:
\begin{lemma}\label{ravel}
Let $U$ be a $\times_A$-Hopf algebra
and assume $A$ is finitely generated
projective as a $U$-module. Then the
cap product with the fundamental class 
$\omega \in H_0(A^*)=A^* \otimes_U A$
defines for all $M \in \umod$ an isomorphism 
$$
		  \cdot \smallfrown \omega :
		  H^0(M) \rightarrow H_0(M
 \otimes
 A^*).
$$
\end{lemma}
\begin{pf}
We have
$\varphi \smallfrown \omega
=\sum_i  (\varphi (1) \otimes_A e^i)
		  \otimes_U e_i$,
and Lemma~\ref{RechtsLinks} identifies
$$
		  H_0(M \otimes A^*)=
		  (M \otimes A^*) \otimes_U A
 \simeq
		  A^*\otimes_U (A \otimes M)
 \simeq A^* \otimes_U M.
$$
In this chain of identifications, $\varphi \smallfrown
 \omega$ is mapped to
\begin{eqnarray}
 \varphi \smallfrown \omega 
&\mapsto&\sum_i  e^i \otimes_U 
		  (e_i \otimes_A \varphi (1))
\mapsto\sum_i  e^i \otimes_U 
		  (e_i\varphi (1))
		  = \sum_i e^i \otimes_U 
		  \varphi (e_i) \nonumber
\end{eqnarray} 
which is identified with 
$\varphi $ under the isomorphism
$\mathrm{Hom}_U(A,M) \simeq 
A^* \otimes_U M$ given by 
$
		  \varphi \mapsto \sum_i e^i
 \otimes_U \varphi (e_i)
$
as in (\ref{sonn}). The claim follows.
\end{pf}
{}\\[-9mm]
  
\subsection{The derived case}
It remains to throw in some 
homological algebra to obtain
Theorem~\ref{main} in general. To shorten 
the presentation, we define:

\begin{dfn}
A module $A$ over a ring $U$ 
is perfect if it admits a finite
resolution by finitely generated
projectives. We call such a module a
duality module if there exists $d \ge 0$
such that $ \mathrm{Ext}^n_U(A,U)=0$
for all $n \neq d$. We abbreviate in
this case $A^*:=\mathrm{Ext}^d_U(A,U)$
and call $d$ the dimension of $A$. 
\end{dfn}

The main remaining step is to prove a
derived version of
Lemma~\ref{cardiff}. One could use a
result of Neeman by which $A \in \umod$
is perfect if and only if 
$\mathrm{Hom}_U(A,\cdot)$
commutes with direct sums 
\cite{Kel:ODGC,Nee:TC}, or 
the Ischebeck spectral
sequence 
which degenerates at $E^2$ if $A$ is a
duality module
\cite{Isch:EDZDFEUT,Kra:PDIHC,Skl:PDARTFEAT}. 
However, we 
include a more elementary and
self-contained proof.

\begin{thm}\label{swansea}
Let $A \in \umod$ 
be a duality module of dimension $d$. 
\begin{enumerate}
\item[1.] The projective dimension 
of $A \in \umod$ is $d$.
\item[2.] $A^*$ is a duality module 
of the same dimension $d$. 
\item[3.] If $P_\bullet \rightarrow A$ is a finitely
generated projective resolution of
length $d$, then 
$P^*_{d-\bullet}=
\mathrm{Hom}_U(P_{d-\bullet},U)$
is a finitely generated projective
resolution of $A^*$ and 
the canonical isomorphism
$$
		 \delta :  M \otimes_U P_i \rightarrow 
		\mathrm{Hom}_U(P^*_i,M),\quad
		m \otimes_U p \mapsto \{\alpha
		\mapsto m \alpha (p)\}
$$
induces for all $U^\mathrm{op}$-modules $M$
a canonical isomorphism
$$
		  \mathrm{Tor}^U_{i}(M,A)
		\rightarrow 
		\mathrm{Ext}^{d-i}_{U^\mathrm{op}}(A^*,M).
$$
\item[4.] There is a canonical isomorphism 
$(A^*)^* \simeq A$.
\end{enumerate} 
\end{thm}
\begin{pf}
Let $P_\bullet \rightarrow A$ be a 
finitely generated projective 
resolution of finite length 
$m \ge 0$ (which exists since 
$A$ is perfect). 
Then the (co)homology of  
$$
		  0 \rightarrow P_0^* \rightarrow \ldots
		  \rightarrow P_m^* \rightarrow 0,\quad
		  P_n^*=\mathrm{Hom}_U(P_n,U)
$$
is 
$ \mathrm{Ext}_U^\bullet(A,U)$, so by 
assumption we have $m \ge d$
and the above complex is
exact except at $P^*_d$ where the homology is 
$A^*=\mathrm{Ext}_U^d(A,U)$. Furthermore, 
all the $P_n^*$ are finitely generated 
projective since the $P_n$ are so 
(Lemma~\ref{cardiff}).

Let $\pi_i$ be the map 
$P^*_i \rightarrow P^*_{i+1}$
and put $K:=\mathrm{ker}\, \pi_{d+1}$. 
By construction, 
\begin{equation}\label{auflnull}
0 \rightarrow K \rightarrow 
P^*_{d+1} \rightarrow \ldots
 \rightarrow P^*_m \rightarrow 0
\end{equation} 
is exact. If one compares this 
exact sequence with the sequence  
$$
		  \ldots \rightarrow 0 \rightarrow
 0 \rightarrow P^*_m \rightarrow
 P^*_m \rightarrow 0
$$
using Schanuel's lemma 
(see \cite[7.1.2]{McCRob:NcNR}), one obtains that
$K$ is projective.

The exactness of $P_\bullet^*$ at 
$P^*_{d+1}$ gives 
$K=\mathrm{im}\, \pi_d$, and 
by the projectivity of $K$, the map 
$\pi_d : P^*_d \rightarrow K \subset 
P^*_{d+1}$ splits so that 
$P^*_d \simeq K \oplus K_\perp$,
$K_\perp:=\mathrm{ker}\, \pi_d$. 
In particular, both 
$K$ and $K_\perp$ are
finitely generated.

It follows from all this that the complex
\begin{equation}\label{auflastern}
		  0 \rightarrow P^*_0 \rightarrow
 \ldots
 \rightarrow P^*_{d-1} \rightarrow
 K_\perp
\rightarrow 0
\end{equation} 
is a finitely generated projective
 resolution of $A^*$: since
$\mathrm{im}\, \pi_{d-1} \subset 
P^*_d$ is contained in 
$\mathrm{ker}\, \pi_d = K_\perp$,
the above complex is still 
exact at $P^*_{d-1}$, and the homology
 at $K_\perp$ is the homology of
 $P^*_\bullet$ at $P^*_d$, that is,
 $A^*$.

Since (\ref{auflnull}) is a finitely
generated projective resolution of $0$
and 
$P^*_{d-\bullet}$ is as a complex a 
direct sum of (\ref{auflastern}) and (a
 shift of) (\ref{auflnull}) we also know
that
 $\mathrm{Ext}^\bullet_{U^\mathrm{op}}
(A^*,M)$ is
for any $M \in \modu$ the (co)homology of
$\mathrm{Hom}_U(P^*_{d-\bullet},M)$.  
By Lemma~\ref{cardiff}, 
this is isomorphic as a chain complex to 
$M \otimes_U P_{d-\bullet}$  
via the isomorphism given in 3., 
and the homology of this
complex is 
$\mathrm{Tor}^U_{d-\bullet}(M,A)$. This
proves 3. The special case 
$M=U$ implies the remaining claims.
\end{pf}

Assume finally
that in the situation of the above
theorem, $U$ is an $A$-biprojective $\times_A$-Hopf
algebra. Since
$P$ is a projective resolution, we have 
$M \otimes_U P \simeq M \lotimes_U P$ 
and $\mathrm{Hom}_{U}(P^*,M) \simeq 
\mathrm{RHom}_U(P^*,M)$, and $\delta$ 
gives an isomorphism between them. 
The fundamental class is defined
to be
$$
		  \omega := \delta^{-1}
		  (\mathrm{id}_{A^*}) \in A^*
		  \lotimes_U A \simeq 
		  P^* \otimes_U A \simeq
		  A^* \otimes_U P,
$$
and Theorem~\ref{swansea} gives immediately: 
\begin{cor}
If $e_1,\ldots,e_n$ and 
$\tilde e^1,\ldots\tilde e^n$
are generators of $A$ and of  
$A^*$, respectively, then
there are 
$e^1,\ldots,e^n \in P^*_0$ and 
$\tilde e_1,\ldots,\tilde e_n \in P_d$ such that 
$$
		  \omega = \sum_i e^i \otimes_U e_i
		  = \sum_i \tilde e^i \otimes_U 
		  \tilde e_i,
$$
and $\delta$ is given by the Yoneda
 product $\cdot \bullet \omega$.
\end{cor}

Theorem~\ref{main} follows now
as in the underived case
(Lemma~\ref{ravel}) working with 
$\mathrm{RHom}_U(A,M)$ and 
$(M \lotimes A^*) \lotimes_U A$ instead of
$H^0(M)=\mathrm{Hom}_U(A,M)$ and 
$H_0(M \otimes A^*)=(M \otimes A^*)
\otimes_U A$: using 
Theorem~\ref{suarez2} and (\ref{weihnacht})
one gets
\begin{eqnarray}
 (M \lotimes A^*) \lotimes_U A 
&\simeq&
		  A^* \lotimes_U (A \lotimes M)
		  \simeq
		  A^* \lotimes_U M \nonumber\\ 
&\simeq&  P^* \lotimes_U M \simeq
		  \mathrm{RHom}_U(P,M) \nonumber\\
&\simeq&
		  \mathrm{RHom}_U(A,M),\nonumber 
\end{eqnarray} 
where we hide the 
reindexing of the complexes for the sake
of better readability (so $P^*$ stands
for
$P^*_{d-\bullet}$, and
$\mathrm{RHom}_U(P,M)$ and
$\mathrm{RHom}_U(A,M)$ are reindexed in
the same way).
This leads to a convergent spectral
sequence
$$
\mathrm{Tor}^U_p(\mathrm{Tor}^A_q(M,A^*),A)
\Rightarrow  
\mathrm{Ext}_U^{d-p-q}(A,M),
$$
and under the last assumption 
of Theorem~\ref{main}
($\mathrm{Tor}^A_q(M,A^*)=0$ for $q>0$)
this spectral sequence degenerates to the claimed
isomorphism.

\providecommand{\bysame}{\leavevmode\hbox to3em{\hrulefill}\thinspace}
\providecommand{\MR}{\relax\ifhmode\unskip\space\fi MR }
\providecommand{\MRhref}[2]{%
  \href{http://www.ams.org/mathscinet-getitem?mr=#1}{#2}
}
\providecommand{\href}[2]{#2}

\end{document}